# Towards 1ULP evaluation of Daubechies Wavelets

NICK THOMPSON, Oak Ridge National Laboratory
JOHN MADDOCK, Boost Software Foundation
GEORGE OSTROUCHOV, JEREMY LOGAN, DAVID PUGMIRE, and SCOTT KLASKY, Oak Ridge National Laboratory

We present algorithms to numerically evaluate Daubechies wavelets and scaling functions to high relative accuracy. These algorithms refine the suggestion of Daubechies and Lagarias to evaluate functions defined by two-scale difference equations using splines; carefully choosing amongst a family of rapidly convergent interpolators which effectively capture all the smoothness present in the function and whose error term admits a small asymptotic constant. We are also able to efficiently compute derivatives, though with a smoothness-induced reduction in accuracy. An implementation is provided in the Boost Software Library.

## 1 INTRODUCTION

The most successful applications of Daubechies wavelets have been in signal and image processing, where the filter coefficients take center stage and the continuous wavelet is usually ignored. In numerical analysis, Daubechies wavelets have been successfully applied in the solution of integral equations and PDE, perhaps most dramatically in construction of a linear scaling density functional theory[1]. In these applications, numerical evaluation of the basis functions is not as important as that of their quadratures, which are particularly simple, as they form an orthonormal basis of $L^2(\mathbb{R})$. When integrating over integral kernels, this simplicity disappears, and much effort has been put into design of bespoke quadrature routines[2] which *avoid* numerical evaluation. Though this is certainly a path to very high performance, we would like to see a lower barrier to entry. A rapidly convergent, widely available, and well-written algorithm for numerical evaluation of these functions affords us the ability to use standard tools like quadrature routines or rendering algorithms generically, while simultaneously reaping all the benefits of the sparsity, such as reduced memory pressure.

The most commonly used numerical software libraries, such as the GSL[25], NAG[29], Arb[26], and Intel MKL do not provide implementations for evaluation of Daubechies wavelets and scaling functions. PyWavelets[21] allows the user to generate a vector of user-provided length whose elements are wavelet and scaling function values evaluated at an equispaced grid, but no interpolation is provided and it is *a priori* unclear how refined the grid should be to achieve a desired accuracy. Maple[27] provides a tool to *plot* wavelets, but no evaluation of the continuous function.

Mathematica[28] supports evaluation of continuous Daubechies wavelets and scaling functions; for example, we may write `WaveletPsi[DaubechiesWavelet[6],x]` to evaluation $_6\psi$. The algorithm linearly interpolates between a small number of equispaced function values, which is not sufficient to recover high accuracy, and derivatives are unavailable.

## 2 MOTIVATION

We began to think about numerical evaluation of Daubechies wavelets while attempting to visualize massive scientific data stored on magnetic tape. If we require lossless storage, but can accept *lossy visualization*, then the multiresolution analysis suggests a way to reorganize our data during the write phase so that only a tiny fraction needs to be read by the renderer. To wit, let $\phi$ be a scaling function and $\psi$ be the associated wavelet, and define

$$\phi_{jk}(x) := \frac{1}{\sqrt{2^j}} \phi(2^{-j} x - k).$$

Since for all $f \in L^2(\mathbb{R})$

$$f = \sum_{j,k \in \mathbb{Z}} \langle f, \psi_{jk} \rangle \psi_{jk} \approx \sum_{k \in \mathbb{Z}} \langle f, \phi_{j_{\max} k} \rangle \phi_{j_{\max}, k} + \sum_{j=j_{\min}}^{j_{\max}} \sum_{k \in \mathbb{Z}} \langle f, \psi_{jk} \rangle \psi_{jk}, \tag{1}$$

we can store the expansion coefficients rather than samples $f(2^{j_{\min}} k)$. By restricting to a finite interval of interest $[a, b]$, we see that $f$ can be approximated at the $j$th level with a fraction $1/2^{j_{\max} - j_{\min} - j}$ of the original data. But for this to work, we must be able to numerically evaluate $\psi_{jk}$; otherwise we have to perform a full inverse transform. The same is true of thresholded representations, where most of the expansion coefficients are set to zero, and hence inverting the transform requires computations on a massive buffer filled primarily with zeros. If we can do all the desired operations directly on the compressed representation of a thresholded transform, we can increase speed and reduce memory pressure.

Of course our motivation does not exhaust the possible uses. For example, if we can accurately evaluate $\psi_{jk}$, then we can evaluate the inner products $\langle f, \psi_{jk} \rangle$ by quadrature. This might be very useful when the data is given on an irregular grid, or when we grow suspicious of the approximation $f(2^j k) \approx 2^{-j/2} \langle f, \phi_{jk} \rangle$ [3] common in the signal processing literature.

## 3 NUMERICAL PRELIMINARIES

Let $\{c_k\}_{k=0}^{2p-1}$ be the filter coefficients of a compactly supported scaling function with $p$ vanishing moments. Evaluation of scaling functions begins with the two-scale difference equation

$$\phi(x) = \sum_{k=0}^{2p-1} c_k \phi(2x - k). \tag{2}$$

This equation, along with the normalization condition

$$\sum_{k \in \mathbb{Z}} \phi(x - k) = 1 \tag{3}$$

*defines* $\phi$. If we evaluate (2) at the integers, we obtain

$$\phi(j) = \sum_{k=0}^{2p-1} c_k \phi(2j - k) = \sum_{k=1}^{2p-2} c_{2j-k} \phi(k), \quad j = 1, 2, \ldots 2p - 2,$$

where the $j$'s are selected from the theorem that $\mathrm{supp}(\phi) = [0, 2p - 1]$ whenever $\mathrm{supp}(\mathbf{c}) = [0, 2p - 1]$[3]. This is an eigenvalue problem $\boldsymbol{\phi} = \mathbf{L}\boldsymbol{\phi}$ with $L_{jk} := c_{2j-k}$[5]. Solving dense eigensystems is in general quite expensive, but this one is easily solved since the value of $p$ of practical interest tends to be on the order of 10. The eigenvector $\boldsymbol{\phi}$ is normalized via (3) which gives us the values of $\phi$ at the integers. Then (2) allows evaluation at half integers as

$$\phi(i/2) = \sum_{k=0}^{2p-1} c_k \phi(i - k).$$

We can continue this recursion to arbitrary depth, obtaining the value of $\phi$ at all *dyadic rationals* $n/2^j$ via

$$\phi(n/2^j) = \sum_{k=0}^{2p-1} c_k \phi(n/2^{j-1} - k).$$

We will call a table of function values given at a set of dyadic rationals a *dyadic grid*. Dyadic grids grow exponentially in $j$, requiring storage of $\sim 2^{j+1}(2p-1)$ floating point numbers. Instead of storing the values at every dyadic rational, we can make a time-space trade-off by using

$$\phi(n/2^j) = \sum_{\ell_1, \ell_2, \ldots, \ell_j} c_{\ell_1} \cdots c_{\ell_m} \phi(n - 2^{j-1}\ell_1 - 2^{j-2}\ell_2 - \cdots - 2\ell_{j-1} - \ell_j)$$

This sum only requires storage of $\phi$ at integers, but then gives a sum which requires $2^{j+1}(2p-1)^2$ terms to evaluate. So in these two extreme cases, we have an algorithm which is exponential in memory but $O(1)$ in evaluation, and another which is exponential in evaluation and $O(1)$ in memory.

We can of course "meet in the middle". We can refine the table up to some $j_{\max}$, and then use (2) once to find

$$\phi\left(\frac{n}{2^{j_{\max}+1}}\right) = \sum_{k=0}^{2p-1} c_k \phi\left(\frac{n}{2^{j_{\max}}} - k\right) \tag{4}$$

or twice to obtain

$$\phi(n/2^{j_{\max}+2}) = \sum_{k_1=0}^{2p-1} c_{k_1} \phi(n/2^{j_{\max}+1} - k_1) = \sum_{k_1,k_2=0}^{2p-1} c_{k_1} c_{k_2} \phi(n/2^{j_{\max}} - 2k_1 - k_2) \tag{5}$$

This method is called "turning the crank", or "cranking". A method very similar in spirit to turning the crank is the cascade algorithm: Given any function $\phi^{(0)}$, the iteration

$$\phi^{(i+1)}(x) = \sum_{k=0}^{2p-1} c_k \phi^{(i)}(2x - k) \tag{6}$$

converges in $L^p$[7]. Unfortunately, neither cascading nor cranking can be expected to produce high relative accuracy. This is because of the nonlocal nature of (4), which gives the value of $\phi$ near a root as a filtered sum of values of $\phi$ far from the root, which are unlikely to be small.

Derivatives can also be exploited to aid numerical evaluation. Differentiating (2) gives

$$\frac{1}{2}\phi'(j) = \sum_{k=0}^{2p-1} c_k \phi'(2j - k), \tag{7}$$

another eigenvalue problem but with corresponding eigenvector $1/2$ rather than 1. (How to normalize the derivatives is described in Oslick et al[8], and explicit values to verify the computation can be found in Lin *et al*[9].) Hence we can generate *two* lookup tables, one for $\phi$, and another for $\phi'$, and use (say)

$$\phi(n/2^j + \epsilon) \approx \phi(n/2^j) + \epsilon \phi'(n/2^j).$$

It is not *a priori* clear that storing a few (hopefully) smaller derivative tables is an efficient use of memory. The derivatives of Daubechies wavelets are strange, as the final derivative is always "weak". For example, $_2\phi$ has a derivative in the sense that the eigenvalue problem $\frac{1}{2} {_2\boldsymbol{\phi}'} = \mathbf{L}_2 \boldsymbol{\phi}'$ has a non-trivial solution, but $_2\phi$ is not continuously differentiable. That said, constructing grids of derivatives is certainly an efficient use of *time*, since they can be computed in parallel.

Despite the bleak prospects of cascading and cranking, a highly refined dyadic grid (along with derivatives, if desired), computed in greater than working precision and evaluated with a local interpolant *can* be expected to produce high relative accuracy. The question is how refined does the grid have to be? We can start by considering the most extreme case: Simply storing values of $\phi$ at *every* floating point representable on the interval $[0, 2p - 1]$. For $_2\phi$, there are roughly a billion 32-bit floats on $[0, 4]$, requiring a lookup table about 4 GB in size. However, a lookup table of $_2\phi$ providing every double precision value at each double-precision representable on $[0, 4]$ occupies well over an exabyte. In order to achieve low relative error with an object occupying a reasonable amount of memory, interpolating the resulting table is required. Daubechies and Lagarias[11] give convergence rates for splines approximating Hölder continuous functions, giving us some intuition about how efficient this approach can be. To wit, if $f$ is $L$ times differentiable with $f^{(L)}$ Hölder continuous with exponent $\alpha$ and $f_j \in C^L$ is a spline which interpolates $f$ at dyadic rationals spaced by $2^{-j}$, then there exists $C > 0$ such that

$$\left\| f - f_j \right\|_\infty \le C \cdot 2^{-(L+\alpha)j}. \tag{8}$$

Jackson's theorem[18] shows that the same rate of convergence is expected for interpolation by trigonometric polynomials–perhaps hinting at the fundamental amount of difficulty in approximating this class of functions. The convergence rate (8), as stated, *looks* exponential, but in fact we should not think about the number of dyadic grid refinements $j$, but the number of elements stored in the table $n = O(2^j)$. Hence the "true" rate of convergence is merely algebraic. The base-2 logarithm of (8) gives a rough estimate of the number of correct bits gained on each dyadic grid refinement: $L + \alpha$. For $_2\phi$, with $L = 0$ and $\alpha \approx 0.55$, this is a very severe constraint. Each doubling in memory can only be expected to deliver roughly *half* a bit of accuracy. Achieving double precision, with 52 mantissa bits, should require 104 dyadic grid refinements–worse than brute force enumeration of every representable! However, for smoother scaling functions, such as $_{15}\phi$, with $L + \alpha \approx 4.5$, we can recover all 52 mantissa bits at $j = 10$, with a lookup table occupying only 238KB. Hence, if we only require 32 or 64 bits of precision on smoother wavelets, the exponential explosion in memory might not be such a pressing issue.

## 4 CHOOSING THE INTERPOLATOR

Because Daubechies wavelets and scaling functions have such unfamiliar smoothness, we did not think we could derive the best interpolator for the dyadic grid from theoretical considerations alone, and instead did a brute force search through every interpolator provided in Boost: piecewise constant, linear, cardinal quadratic, cubic, and quintic $B$-splines[13], first, second, and third-order Taylor expansions, modified Akima interpolation, PCHIP[14], Whittaker-Shannon[16], trigonometric[18], as well as cubic, quintic, and septic Hermite splines. For interpolators which support non-equispaced abscissas, we also tried precomputing a subset of the roots of the scaling function and sprinkled those values in; the idea being that partitioning the support into positive and negative subintervals might allow evaluation with only positive or only negative terms. This fixed some gross sign errors at low refinements, but didn't reduce the sup-norm error significantly. In addition, it made the abscissa spacing irregular, forcing $O(\log(n))$ evaluation as binary search is required to identify the correct subinterval, rather than $O(1)$ evaluation possible with an equispaced table.

Once we had a decent interpolator, we were able to "bootstrap" to the Remez algorithm[17], Chebyshev series[19], and barycentric rational interpolation[20], all of which require that a function can be evaluated at any point, rather than only a small number of dyadic rationals. In retrospect this was akin to expecting a miracle, since these require $O(n)$ evaluation, and hence we would have required a vast increase in convergence rates to make any of them competitive with the $O(1)$ methods. We are therefore relieved to report dramatic failure of all of them: The resulting linear system solved by the Remez algorithm was very poorly conditioned, and we could not push the absolute error below $10^{-2}$. The Chebyshev series took up a massive amount of memory, and the barycentric rational interpolation converged slowly.

For each interpolator, $\log_2\left(\|\phi - \phi_j\|_\infty\right)$ plotted against the dyadic grid refinement $j$ was so perfectly linear that we thought it sensible to compute linear regressions on

$$\log_2\left(\|\phi - \phi_j\|_\infty\right) \approx \log_2(\tilde{C}) - \tilde{\alpha} j$$

where $j$ is the number of dyadic grid refinements, $\tilde{C}$ is the empirical asymptotic constant, $\tilde{\alpha}$ is an empirical convergence rate that is exceedingly close to the actual Hölder exponent $L + \alpha$ if the interpolator is sufficiently smooth. (These regressions are of course suggested by (8).) Since so many of the interpolators are smooth enough to capture the smoothness of the wavelets and scaling functions for the range of vanishing moments $p$ we considered, the real problem is to find the interpolator which admits the smallest asymptotic constant $\tilde{C}$. On this count, simple Hermite splines are a clear winner, with Schoenberg's $B$-splines a distant second. (For those unfamiliar with Hermite splines, or simply confused by which of the many definitions of "Hermite spline" is being used here, we have added a short description to the appendix.) Hermite splines require not only the dyadic grid, but also as certain number of derivatives, so there are two considerations: Convergence as $j \to \infty$ and memory usage. Fortunately, no difficult decisions need to be made; the convergence of the Hermite splines is so dramatically faster than the interpolators requiring only the values of the dyadic grid that total memory consumption was less. In addition, the derivatives of the Hermite splines are vastly more accurate than the derivatives of the interpolants which do not use derivative data, so applications which require accurate derivatives are not forced to choose between memory and accuracy.

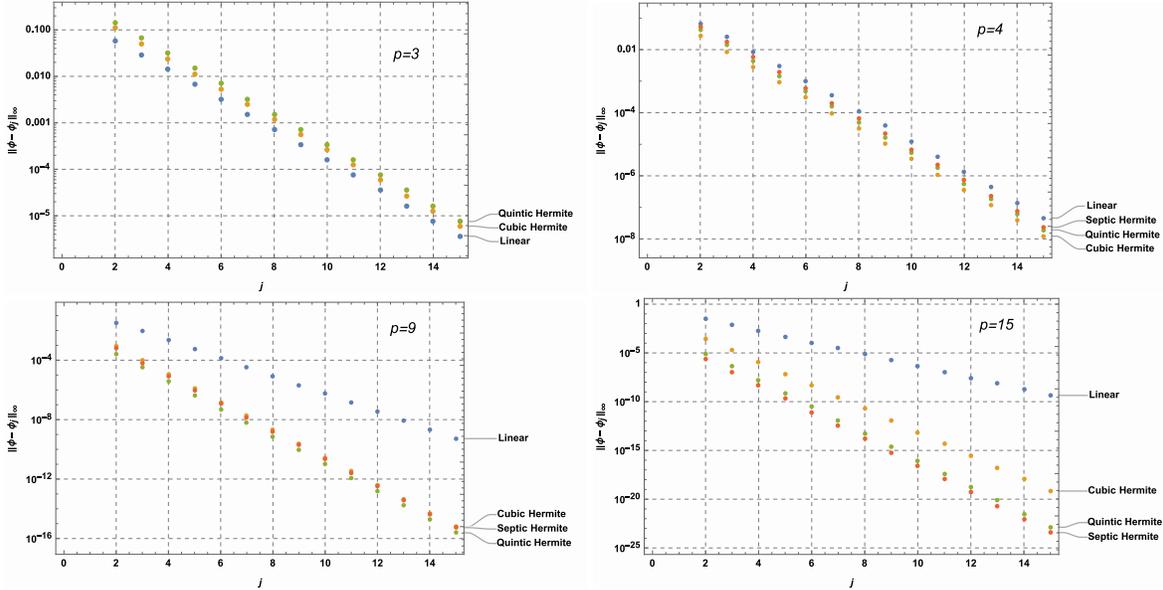

Fig. 1. Different interpolators converge faster for scaling functions with a different number of vanishing moments. For the best interpolator, the slopes of the lines (in $\log_2$) are very nearly the Hölder exponents.

## 5 THE CHALLENGE OF TWO VANISHING MOMENTS

The Daubechies scaling function $_2\phi$, with the least smoothness, is the most difficult to evaluate. We found amongst standard interpolators that linear interpolation of the $_2\phi$ dyadic grid was the best, but nevertheless an enormous amount of memory was required to produce acceptable accuracy. Our implementation therefore exploits a famous property of ${}_2\phi$[10]: that it is right-differentiable at dyadic rationals, but not left-differentiable, with the differentiability failing as a root singularity with exponent $\alpha := 2 - \ln(1+\sqrt{3})/\ln(2) \approx 0.55$[12]. If $x \in [i/2^{j_{\max}}, (i+1)/2^{j_{\max}}] =: [x_i, x_{i+1}]$, then the interpolation problem is to find $f$ such that

$$\lim_{x \nearrow x_{i+1}} f(x) = {}_2\phi(x_{i+1}), \quad \lim_{x \nearrow x_{i+1}} f'(x) = {}_2\phi'(x_{i+1}), \quad \lim_{x \searrow x_i} f(x) = {}_2\phi(x_i), \quad \lim_{x \searrow x_i} f(x) \sim C(x-x_i)^\alpha + {}_2\phi(x_i)$$

| $p$ | Best interpolator | $\log_2(\lVert {}_p\phi - {}_p\phi_j \rVert_\infty)$ | 64 bit evaluation (ns) | RAM (relative, Mb) | RAM (absolute, kb) |
|---|---|---|---|---|---|
| 2 | Matched Hölder | -2.15 - 0.55j | 2.6 | 50 | 6292 |
| 3 | Linear | -1.85 - 1.08j | 1.7 | 84 | 20972 |
| 4 | Cubic Hermite | -1.98 - 1.62j | 3.5 | 117 | 3670 |
| 5 | Cubic Hermite | -3.93 - 1.98j | 3.5 | 9.4 | 148 |
| 6 | Quintic Hermite | -4.96 - 2.20j | 5.5 | 8.6 | 34 |
| 7 | Quintic Hermite | -4.76 - 2.46j | 5.5 | 5.1 | 20 |
| 8 | Quintic Hermite | -4.99 - 2.77j | 5.5 | 2.9 | 23 |
| 9 | Quintic Hermite | -5.72 - 3.08j | 5.5 | 1.6 | 26 |
| 10 | Septic Hermite | -7.21 - 3.36j | 8.4 | 1.2 | 9.7 |
| 11 | Septic Hermite | -7.89 - 3.61j | 8.4 | .68 | 10.8 |
| 12 | Septic Hermite | -8.66 - 3.86j | 8.4 | .75 | 5.9 |
| 13 | Septic Hermite | -9.29 - 4.10j | 8.4 | .82 | 6.4 |
| 14 | Septic Hermite | -9.14 - 4.32j | 8.4 | .88 | 6.9 |
| 15 | Septic Hermite | -9.47 - 4.56j | 8.4 | .95 | 7.5 |
| 16 | Septic Hermite | -10.1 - 4.80j | 8.4 | 1.01 | 4.0 |
| 17 | Septic Hermite | -10.8 - 5.02j | 8.4 | 1.08 | 4.2 |
| 18 | Septic Hermite | -10.6 - 5.24j | 8.4 | 1.14 | 4.5 |
| 19 | Septic Hermite | -10.9 - 5.46j | 8.4 | 1.2 | 4.8 |

Table 1. Best interpolators for each ${}_p\phi$. The optimal Hölder exponents of ${}_2\phi$ ($\alpha = 0.5500$), ${}_3\phi$ ($\alpha = 1.08$), and ${}_4\phi$ ($\alpha = 1.619$) from Daubechies and Lagarias[12] show up clearly in the empirical convergence, showing that our chosen interpolators are capturing all the smoothness that exists in the function. Timings were measured using Google Benchmark on an Intel Core i7-7820X CPU, compiled with clang 9.0.0 using `-ffast-math -O3 -march=native` using a favorable access pattern (increasing abscissas). More realistic access patterns are considered in the examples.

Using the ansätz

$$f(x) = {}_2\phi(x_i) + c_1 \left(\frac{x - x_i}{h}\right) + c_2 \left(\frac{x - x_i}{h}\right)^\alpha$$

we obtain very simple expressions for $c_1$ and $c_2$:

$$c_1 = \frac{1}{1-\alpha}\left[{}_2\phi'(x_{i+1})h - \alpha({}_2\phi(x_{i+1}) - {}_2\phi(x_i))\right]$$

$$c_2 = \frac{1}{1-\alpha}\left[-{}_2\phi'(x_{i+1})h + \phi(x_{i+1}) - \phi(x_i)\right],$$

where $h := x_{i+1} - x_i = 1/2^{j_{\max}}$. We refer to this scheme as "matched Hölder interpolation". If we use a somewhat smaller value of $\alpha$ than the true exponent of 0.5500, the values of the function are more accurate over the interpolation interval, while capturing the asymptote only slightly less well. Since calls to `std::pow` are expensive, we pushed the exponent down to precisely 0.5, allowing a call to `std::sqrt`, which is an order of magnitude faster. This scheme reduces the sup norm error by a factor of 3 relative to linear interpolation; a very welcome improvement saving two grid refinements at a given accuracy.

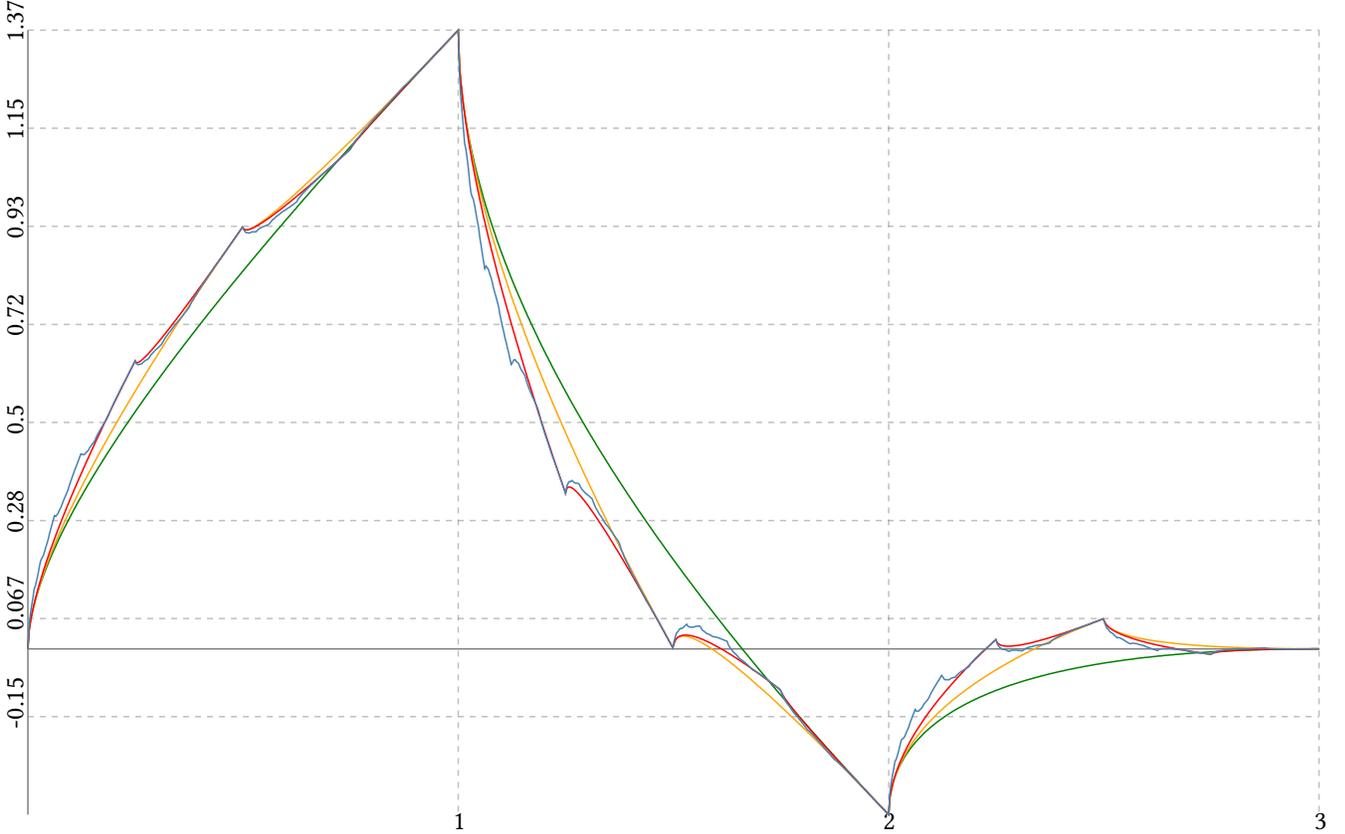

Fig. 2. $_2\phi$ computed at 0 (green), 1 (orange), 2 (red), and 24 (blue) dyadic grid refinements, and interpolated with matched Hölder interpolation. The improvement over linear interpolation of the kinks is obvious.

## 6 CHOOSING THE REFINEMENT

With the fastest converging interpolator in hand, we now determine how refined the dyadic grid must be in order to produce accurate evaluation. The main surprise is how much more difficult controlling relative error is than controlling absolute error. For example, if we wish to evaluate $_{15}\phi$ to an absolute error on the order of the 32-bit float unit roundoff, then we use Table 1 to find $j$ such that $\left\Vert \phi - \phi_j \right\Vert_\infty \approx 2^{-9.47-4.56j} \leq 2^{-23}$, giving $j \geq 3$. However, if we wish have reasonable control over the ULP distance away from the roots, we must take $j = 14$, requiring 2000 times more memory! This is entirely the fault of the condition number of function evaluation

$$\mathrm{cond}[f](x) := \left| \frac{xf'(x)}{f(x)} \right|,$$

which grows superexponentially for Daubechies scaling functions towards the right end of the support (see Figure 3). The whole enterprise seems more natural in *fixed point*, where the condition number is simply $|f'(x)|$, which is bounded for all wavelets and scaling functions except $_2\phi$. This means we can achieve low absolute error with much greater ease than low relative error.

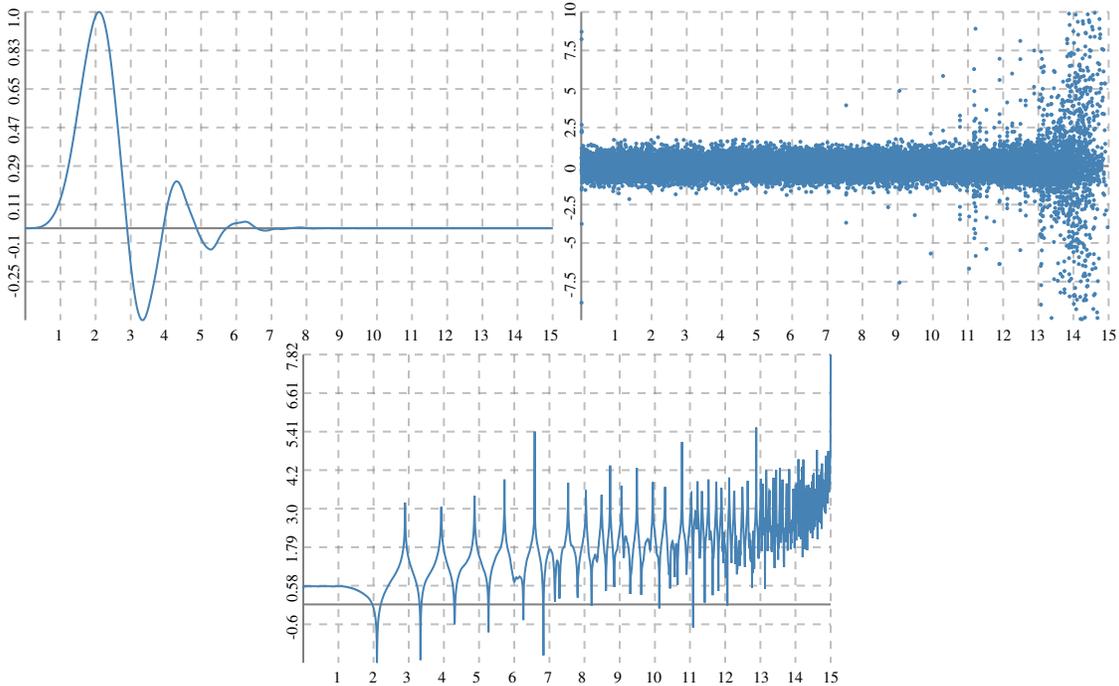

Fig. 3. $_8\phi$ (top left), $\log_{10}(\mathrm{cond}(_8\phi))$ (bottom left), and an ULPs plot of the double precision Boost implementation at 21 refinements. Superexponential growth of the condition number towards the right end of the support is clear, with an obvious effect on the ULPs plot.

We decided on the following criteria for choosing the default number of dyadic grid refinements: If the error to be expected from rounding the argument to the nearest floating point representable is on the same order as the approximation error of the interpolator, then we've done all we need to do. This translates into

$$\left| \frac{_p\phi(x) - _p\phi_j(x)}{_p\phi(x)} \right| \lesssim \mu \, \mathrm{cond}[_p\phi](x) \quad \forall x \in \mathrm{supp}(\phi)$$

whenever $_p\phi(x) \neq 0$. Determining a more rigorous criteria for selection of the default refinements is very challenging. Each datum in the dyadic grid is contaminated by error of order $\mu_P$, where $\mu_P$ is the unit roundoff of the precise floating point type, before casting. But the dense set of roots towards the right side of the support means than many values of the scaling function are denormalized in the target precision; to get these correct would require a refined floating point more than four times the bit budget of the target. Our experience casting from four times the target precision was poor; we gained the intellectual satisfaction of accuracy when $_p\phi(x)$ evaluated to a denorm, at the cost of *extreme* sluggishness constructing the dyadic grid. In addition, the probability of hitting such an abscissa is very low; ULPs plots[6] look very reasonable using 80-bit long double to create a double precision dyadic grid. Eventually, we decided to let the ULPs plot determine the default grid refinements; these keep the evaluation to ~1.5 ULPs in the well-conditioned regions. However, the number of grid refinements is exposed to the user, and the ability to create ULPs plots from them is also provided. We also provide an default number of grid refinements for applications which only require control of absolute error; this is common in computer graphics since pixel spacing is uniform.

## 7 THE WAVELET

Once we have the ability to rapidly evaluate the scaling function, the wavelet can be evaluated via

$$\psi(x) = \sum_{k=0}^{2p-1} (-1)^{k+1} c_k \phi(2x + k - 1).$$

There are a few reasons why we shouldn't apply this formula directly: First, (1) shows that we expect the wavelet to be evaluated with much greater frequency than the scaling function. So if $\phi$ evaluates in 10 nanoseconds, then $\psi$ evaluates in roughly $20p$ nanoseconds, which for commonly-used values of $p$ is unacceptably costly. Second, the sum alternates, so we cannot expect to achieve high relative accuracy. We are therefore forced into the same situation we were in before: Precomputation of a dyadic grid in higher precision, and then interpolating. However, since the sum is finite, $_p\phi$ and $_p\psi$ have identical smoothness, and hence we can use the same interpolator for $_p\phi$ as for $_p\psi$.

## 8 SOFTWARE IMPLEMENTATION

The documentation of our implementation is provided at boost.org, but for completeness we give a short description of our API here. The simplest use generates a C++ class with the default constructor, picking the default number of grid refinements as discussed above. The constructor computes the dyadic grid and its derivatives in parallel, but is still expensive. Hence we have used the PIMPL pattern to allow the class to be shared between threads without copying the underlying dyadic grid data. We define the call operator on this class to give the feel of calling a standard special function, and use `.prime` and `.double_prime` to evaluate derivatives; these calls are threadsafe. The number of vanishing moments is a template parameter, and hence if we call `psi.double_prime(x)` on a wavelet which doesn't have a second derivative, we get an informative compile-time error.

We made much effort to avoid using a class and instead deploy a free function, but were unsuccessful. The obvious way to implement a free function is to precompute the dyadic grids required for each $p$ using extended precision arithmetic and deploy them in a `.cpp` file. Unfortunately, while modern compilers proved perfectly capable of rapidly parsing the large source files so generated, we estimated they would occupy well in excess of 100Mb for *every $p$* we wished to deploy. Adding files of such size to an open source repository was deemed anti-social at best, and the approach was abandoned. However, we provide a routine `boost::math::daubechies_scaling_dyadic_grid` which allows users to precompute the dyadic grids up to any refinement, any derivative, and up to quad precision. This would allow a porting to the C language, and facilitate writing implementations which do not have the same requirements for source distribution to achieve trivial evaluation.

## 9 EXAMPLES

We now verify that our implementation is indeed usable as a primitive for other generic numerical algorithms by computing the continuous wavelet transform with respect to $_p\psi$

$$\mathcal{W}_p[f](s,t) := \frac{1}{\sqrt{|s|}} \int_{-\infty}^{\infty} f(x) \psi_p\left(\frac{x-t}{s}\right) dx = \sqrt{|s|} \int_{-p+1}^{p} f(su+t) \psi_p(u) du.$$

The latter expression, along with the Euler-MacLaurin summation formula suggests using trapezoidal quadrature, which has the benefit of being arbitrarily adaptive, unlike (say) Gauss-Kronrod which admits only a single

```cpp
#include <boost/math/special_functions/daubechies_scaling.hpp>
#include <boost/math/special_functions/daubechies_wavelet.hpp>
using boost::math::daubechies_scaling;
using boost::math::daubechies_wavelet;
int main() {
    constexpr const int p = 8;
    auto phi = daubechies_scaling<double, p>();

    std::cout << "phi(0.25) = " << phi(0.25) << "\n";
    std::cout << "phi'(0.25) = " << phi.prime(0.25) << "\n";
    std::cout << "phi''(0.25) = " << phi.double_prime(0.25) << "\n";
    auto [a, b] = phi.support();
    std::cout << "Support: [" << a << ", " << b << "]\n";
    // The wavelet:
    auto psi = daubechies_wavelet<double, p>();
    // ...
}
```

Fig. 4. An example usage of the Daubechies wavelet and scaling functions provided by Boost.Math. Despite the precomputation of the dyadic grid forcing us to use a class rather than a free function, a class allows us to provide additional metadata beyond numerical evaluation.

adaptive step. Although the purpose is to demonstrate the ease and efficiency with which we can now perform a useful task, we have encapsulated this knowledge into a quadrature routine, whose use is demonstrated in Figure 5. A scalogram computed from the computation is displayed in Figure 6.

```cpp
#include <boost/math/quadrature/wavelet_transforms.hpp>
using boost::math::quadrature::daubechies_wavelet_transform;
int main() {
    auto f = [](double x) {
        if(t==0) {
            return double(0);
        }
        return std::sin(1.3/x);
    };
    auto Wf = daubechies_wavelet_transform<decltype(f), double, 8>(f);
    std::cout << Wf(/*s = */ 2.0, /*t = */ 1.5) << "\n";
}
```

Fig. 5. Computing the continuous wavelet transform of $\sin(1.3/x)$ with respect to $_8\psi$.

Our final example returns to the original motivation: Evaluation of the wavelet series (1) to create graphs. This serves as a minimal working example of *sparse computer graphics*, which assumes that decompression out of a sparse basis into traditional pointsets will often be impossible due to memory constraints. Therefore creating a graphics pipeline which does all requisite operations directly on the compressed representation is essential. In order to demonstrate that our numerical evaluation is useful for this task, we use the "bumps" function defined by Donoho[30] and compute its expansion coefficients $\langle f, {}_p\phi_{jk}\rangle$ and $\langle f, {}_p\psi_{jk}\rangle$ at scales $j = -11, \ldots, 0$. We choose

```cpp
template<class Real>
Real bumps(Real x) {
    if (x <= 0 || x >= 1) {
        return 0;
    }
    static constexpr std::array<Real, 11> t{.1, .13, .15, .23, .25, .40, .44, .65, .76, .78, .81};
    static constexpr std::array<Real, 11> h{4.0, 5.0, 3.0, 4.0, 5.0, 4.2, 2.1, 4.3, 3.1, 5.1, 4.2};
    static constexpr std::array<Real, 11> w{.005, .005, .006, .01, .01, .03, .01, .01, .005, .008, .005};

    Real f = 0;
    for (size_t i = 0; i < 11; ++i) {
        Real z = abs((x-t[i])/w[i]);
        f += h[i]*pow(1+z, -4);
    }
    return f;
}
```

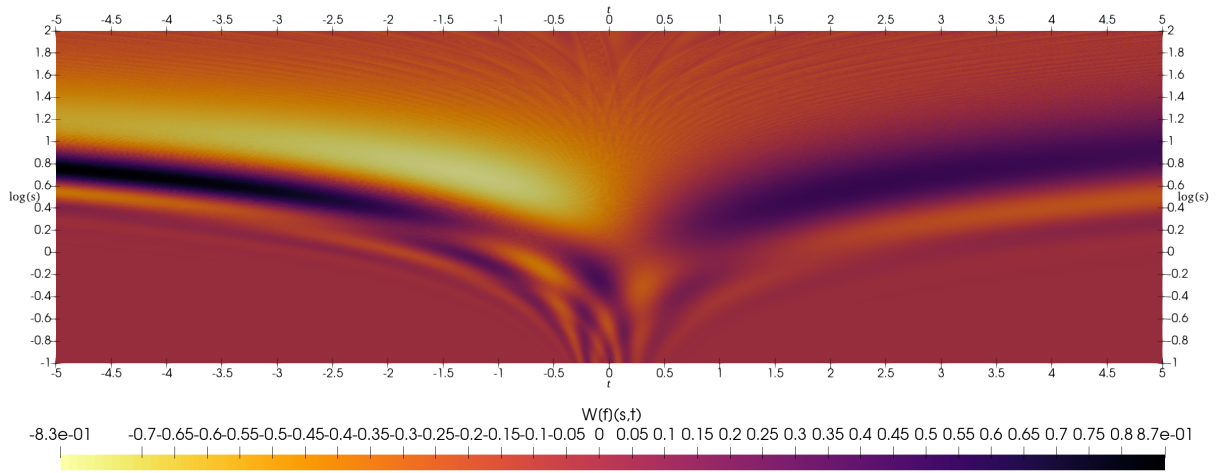

Fig. 6. Wavelet transform of $\sin(1/t)$ computed with respect to $_8\psi$. The pixels are calculated via trapezoidal quadrature and the image is rendered using VTK-m[22] and Paraview.

Fig. 7. A C++ implementation of Donoho's bumps function.

the number of vanishing moment $p$ to maximize the Hoyer sparsity[31] of the expansion, and threshold at a level $10^{-3}$. The results are displayed in Figure 8. Unlike a Fourier or Chebyshev series, we cannot use an evaluation of $_p\psi$ to seed a recurrence that might allow a rapid evaluation of the entire series. However, we have something better: Compact support. This allows us to only read the coefficients whose associated basis function has support that includes the requested abscissa. Using this property of the wavelet series, we can construct the graph at 850,000 pixels/second single threaded, not much worse than using the original definition of Donoho's bumps function, which evaluates at 2.5 million pixels/second.

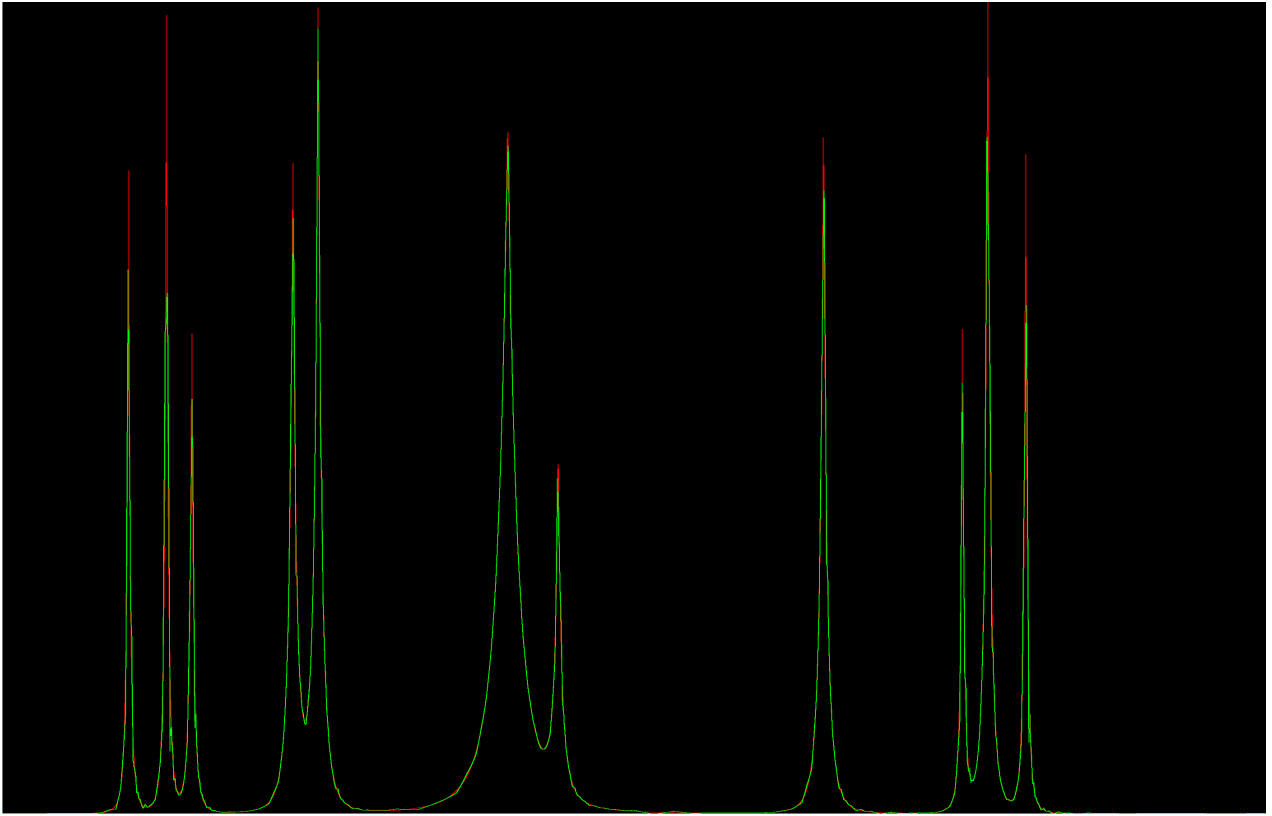

Fig. 8. The bumps function of Donoho (red) and the thresholded wavelet expansion computed with respect to $p = 3$ vanishing moment Daubechies wavelet and scaling function, at scales $j = -8, \ldots, 0$. (The pixels have neither been aliased nor thickened in order to more precisely measure the "fundamental" amount of computation required to plot these functions.) The threshold is chosen as $10^{-3}$. The pre-thresholded, pixel identical wavelet series requires 11 scales and occupies 33kb, the thresholded, 8 scale representation occupies 1.4kb.

In addition to sparsity, the multiresolution analysis affords us the ability to do *progressive rendering*, where we render as many scales as possible in a given timelimit. If $f$ is $C^n$ with $f^{(n)}$ Hölder continuous with exponent $\alpha$, then the expansion coefficients $\langle f, {}_p\psi_{jk}\rangle$ decay as

$$\left|\langle f, {}_p\psi_{jk}\rangle\right| \leq C \cdot 2^{(n+1/2+\alpha)j}$$

as $j \to -\infty$ [4]. This suggests that for slow network transfer or magnetic tape, we should progressively render starting at $j_{\max}$ and repaint as scales are delivered.

## 10 CONCLUSIONS

We have achieved fast and accurate evaluation of the Daubechies wavelets and scaling functions. However, our implementation does not satisfy the same standards as the other special functions in math.h. We have not achieved trivial evaluation with a free function: Our algorithm requires that we must precompute the dyadic grid in the constructor of the `daubechies_scaling` class, in higher precision than the target precision. Second, the algebraic convergence in memory puts 1ULP evaluation in 80 bit `long double` and 128 bit quad precision out of reach of modern consumer hardware. We require either completely new ideas, or hardness results about approximating Daubechies wavelets and scaling functions. We suspect hardness, and believe proof of such results may be within reach of current tools. For example, Alt [23] has demonstrated that the complexity of evaluation of the elementary transcendental functions is equivalent to division with respect to Boolean circuit depth, with $j$ correct bits taking $O(\log(j))$ space. Our empirical results indicate recovering $j$ correct bits for Daubechies wavelets requires $O(2^j)$ space. We hope that our examples have convinced the reader that the fundamental complexity of numerical evaluation of Daubechies wavelets is in fact a problem of applied interest, and those with the correct background can either prove or disprove our conjecture.

The principle success is speed of evaluation. Our algorithm delivers values in roughly the same time it takes to evaluate `std::sin`, and faster than evaluation of `std::pow`. In addition, all derivatives are also available on the same timescale, with a loss of accuracy commensurate with the decrease in smoothness.

*All figures, benchmarks, and computational experiments in this paper can be regenerated by running the code here.*

14 • Nick Thompson, John Maddock, George Ostrouchov, Jeremy Logan, David Pugmire, and Scott Klasky

## A  HERMITE SPLINES

For less mathematically inclined programmers who might be tasked with replicating these algorithms in other programming languages, we give explicit expressions for the Hermite splines, which were found to be the most effective interpolators for the dyadic grids.

Suppose we are given a list of abscissas $\mathbf{x} \in \mathbb{R}^n$, a list of ordinates $\mathbf{y} \in \mathbb{R}^n$, and a list of derivatives $\mathbf{v} \in \mathbb{R}^n$, and wish to construct function $P \in C^1[x_0, x_{n-1}]$ which interpolates them. On each subinterval $[x_i, x_{i+1}]$, we have four constraints, which we fit with a cubic polynomial. We write this as

$$P(x_i) = y_i, \quad P'(x_i) = v_i, \quad P(x_{i+1}) = y_{i+1}, \quad P'(x_{i+1}) = v_{i+1}$$

Let $\Delta x := x_{i+1} - x_i$ and $t(x) := (x - x_i)/\Delta x$. We can restate the problem on the interval $[0, 1]$ via

$$p_i(0) = y_i, \quad p_i'(0) = v_i \Delta x, \quad p_i(1) = y_{i+1}, \quad p_i'(1) = v_{i+1} \Delta x$$

where $P(x) = p_i\left(\frac{x-x_i}{\Delta x}\right)$. We now write

$$p_i(t) = h_{00}(t) y_i + h_{10}(t) v_i \Delta x + h_{01}(t) y_{i+1} + h_{11}(t) v_{i+1} \Delta x$$

where $h_{ij}$ are each cubic polynomials. The constraints on the interpolator $P$ can be translated into constraints on the $h_{ij}$ via

$$\begin{aligned}
h_{00}(0) &= 1, & h_{00}'(0) &= 0, & h_{00}(1) &= 0, & h_{00}'(1) &= 0 \\
h_{10}(0) &= 0, & h_{10}'(0) &= 1, & h_{10}(1) &= 0, & h_{10}'(1) &= 0 \\
h_{01}(0) &= 0, & h_{01}'(0) &= 0, & h_{01}(1) &= 1, & h_{01}'(1) &= 0 \\
h_{11}(0) &= 0, & h_{11}'(0) &= 0, & h_{11}(1) &= 0, & h_{11}'(1) &= 1
\end{aligned}$$

The reader may verify that these constraints are satisfied by

$$h_{00}(t) = 2t^3 - 3t^2 + 1, \, h_{10}(t) = t^3 - 2t^2 + t, \, h_{01}(t) = -2t^3 + 3t^2, \, h_{11}(t) = t^3 - t^2$$

and for completeness

$$p_i(t) = (2t^3 - 3t^2 + 1)y_i + t(t^2 - 2t + 1)v_i \Delta x + t^2(3 - 2t)y_{i+1} + t^2(t - 1)v_{i+1} \Delta x$$

The first derivative of the cubic Hermite spline is available and useful. For $x \in [x_i, x_{i+1})$ we have

$$P'(x) = p_i'(t)/\Delta x = \frac{6t(1 - t)}{\Delta x}(y_{i+1} - y_i) + (3t^2 - 4t + 1)v_i + t(3t - 2)v_{i+1}$$

The complexity of evaluating this spline is $O(\log(n))$ if the abscissas **x** are irregularly spaced (since we must bisect into the correct subinterval $[x_i, x_{i+1}]$), and constant time if they are equally spaced. Note that since every evaluation fetches $y_i, v_i, y_{i+1}, v_{i+1}$, then we should use the "array of structs" data layout rather than the "struct of arrays", i.e., we should lay our data out as

$$\{(y_0, v_0), (y_1, v_1), \ldots, (y_{n-1}, v_{n-1})\}$$

rather than

$$\{y_0, y_1, \ldots, y_{n-1}\}, \{v_0, v_1, \ldots, v_{n-1}\}.$$

It is more effort to find explicit expressions for the quintic and septic Hermite spline basis functions. Nothing about these expressions challenges the imagination, but it is tedious to work them out, and for that reason we write them down here.

We begin with the quintic Hermite spline. This interpolator requires positions **y**, velocities **v**, and accelerations **a** defined on set of abscissas **x**. We wish to interpolate this data with a function $P \in C^2[x_0, x_{n-1}]$. On each subinterval $[x_i, x_{i+1}) = [x_i, x_i + \Delta x)$ we need to match 6 conditions, giving rise to a quintic polynomial. Again we write the restriction of $P$ to the subinterval as $P(x) = p_i((x - x_i)/\Delta x) = p(t)$, giving the constraints

$$p_i(0) = y_i, \quad p_i'(0) = v_i \Delta x, \quad p_i''(0) = a_i \Delta x^2, \quad p_i(1) = y_{i+1}, \quad p_i'(1) = v_{i+1} \Delta x, \quad p_i''(1) = a_{i+1} \Delta x^2$$

We write

$$p_i(t) = h_{00}(t) y_i + h_{10}(t) v_i \Delta x + h_{20}(t) a_i \Delta x^2/2 + h_{01}(t) y_{i+1} + h_{11}(t) v_{i+1} \Delta x + h_{21}(t) a_{i+1} \Delta x^2/2$$

whereupon

$$\begin{aligned}
h_{00}(0) &= 1, & h_{00}'(0) &= 0, & h_{00}''(0) &= 0, & h_{00}(1) &= 0, & h_{00}'(1) &= 0, & h_{00}''(1) &= 0 \\
h_{10}(0) &= 0, & h_{10}'(0) &= 1, & h_{10}''(0) &= 0, & h_{10}(1) &= 0, & h_{10}'(1) &= 0, & h_{10}''(1) &= 0 \\
h_{20}(0) &= 0, & h_{20}'(0) &= 0, & h_{20}''(0) &= 2, & h_{20}(1) &= 0, & h_{20}'(1) &= 0, & h_{20}''(1) &= 0 \\
h_{01}(0) &= 0, & h_{01}'(0) &= 0, & h_{01}''(0) &= 0, & h_{01}(1) &= 1, & h_{01}'(1) &= 0, & h_{01}''(1) &= 0 \\
h_{11}(0) &= 0, & h_{11}'(0) &= 0, & h_{11}''(0) &= 0, & h_{11}(1) &= 0, & h_{11}'(1) &= 1, & h_{11}''(1) &= 0 \\
h_{21}(0) &= 0, & h_{21}'(0) &= 0, & h_{21}''(0) &= 0, & h_{21}(1) &= 0, & h_{21}'(1) &= 0, & h_{21}''(1) &= 2
\end{aligned}$$

The reader may verify that these equations are satisfied by

$$\begin{aligned}
h_{00}(t) &= 1 - 10t^3 + 15t^4 - 6t^5 \\
h_{10}(t) &= t - 6t^3 + 8t^4 - 3t^5 \\
h_{20}(t) &= t^2 - 3t^3 + 3t^4 - t^5 \\
h_{01}(t) &= 10t^3 - 15t^4 + 6t^5 \\
h_{11}(t) &= -4t^3 + 7t^4 - 3t^5 \\
h_{21}(t) &= t^3 - 2t^4 + t^5
\end{aligned}$$

The derivative on $[x_i, x_{i+1})$ is

$$P'(x) = (30t^2 - 60t^3 + 30t^4)\frac{(y_{i+1} - y_i)}{\Delta x} + (1 - 18t^2 + 32t^3 - 15t^4)v_i - (12t^2 - 28t^3 + 15t^4)v_{i+1}$$
$$+ \frac{\Delta x}{2}\left[(2t - 9t^2 + 12t^3 - 5t^4)a_i + (3t^2 - 8t^3 + 5t^4)a_{i+1}\right]$$

and the second derivative is

$$P''(x) = (60t - 180t^2 + 120t^3)(y_{i+1} - y_i)/\Delta x^2 + (-36t + 96t^2 - 60t^3)v_i/\Delta x - (24t - 84t^2 + 60t^3)v_{i+1}/\Delta x$$
$$+ (1 - 9t + 18t^2 - 10t^3)a_i + (3t - 12t^2 + 10t^3)a_{i+1}$$

Our final interpolator is the septic. This requires the third derivative of position (jerk), which we denote by **j**. Using the same notation as before,

$$p_i(t) = h_{00}(t)y_i + h_{10}(t)v_i\Delta x + h_{20}(t)a_i\Delta x^2/2 + h_{30}(t)j_i\Delta x^3/6$$
$$+ h_{01}(t)y_{i+1} + h_{11}(t)v_{i+1}\Delta x + h_{21}(t)a_{i+1}\Delta x^2/2 + h_{31}(t)j_{i+1}\Delta x^3/6$$

whereupon

$$h_{00}(0) = 1, \quad h'_{00}(0) = 0, \quad h''_{00}(0) = 0, \quad h^{(3)}_{00}(0) = 0, \quad h_{00}(1) = 0, \quad h'_{00}(1) = 0, \quad h''_{00}(1) = 0, \quad h^{(3)}_{00}(1) = 0$$
$$h_{10}(0) = 0, \quad h'_{10}(0) = 1, \quad h''_{10}(0) = 0, \quad h^{(3)}_{10}(0) = 0, \quad h_{10}(1) = 0, \quad h'_{10}(1) = 0, \quad h''_{10}(1) = 0, \quad h^{(3)}_{10}(1) = 0$$
$$h_{20}(0) = 0, \quad h'_{20}(0) = 0, \quad h''_{20}(0) = 2, \quad h^{(3)}_{20}(0) = 0, \quad h_{20}(1) = 0, \quad h'_{20}(1) = 0, \quad h''_{20}(1) = 0, \quad h^{(3)}_{20}(1) = 0$$
$$h_{30}(0) = 0, \quad h'_{30}(0) = 0, \quad h''_{30}(0) = 0, \quad h^{(3)}_{30}(0) = 6, \quad h_{30}(1) = 0, \quad h'_{30}(1) = 0, \quad h''_{30}(1) = 0, \quad h^{(3)}_{30}(1) = 0$$
$$h_{01}(0) = 0, \quad h'_{01}(0) = 0, \quad h''_{01}(0) = 0, \quad h^{(3)}_{01}(0) = 0, \quad h_{01}(1) = 1, \quad h'_{01}(1) = 0, \quad h''_{01}(1) = 0, \quad h^{(3)}_{01}(1) = 0$$
$$h_{11}(0) = 0, \quad h'_{11}(0) = 0, \quad h''_{11}(0) = 0, \quad h^{(3)}_{11}(0) = 0, \quad h_{11}(1) = 0, \quad h'_{11}(1) = 1, \quad h''_{11}(1) = 0, \quad h^{(3)}_{11}(1) = 0$$
$$h_{21}(0) = 0, \quad h'_{21}(0) = 0, \quad h''_{21}(0) = 0, \quad h^{(3)}_{21}(0) = 0, \quad h_{21}(1) = 0, \quad h'_{21}(1) = 0, \quad h''_{21}(1) = 2, \quad h^{(3)}_{21}(1) = 0$$
$$h_{31}(0) = 0, \quad h'_{31}(0) = 0, \quad h''_{31}(0) = 0, \quad h^{(3)}_{31}(0) = 0, \quad h_{31}(1) = 0, \quad h'_{31}(1) = 0, \quad h''_{31}(1) = 0, \quad h^{(3)}_{31}(1) = 6$$

The solution of these equations are

$$h_{00}(t) = 1 - 35t^4 + 84t^5 - 70t^6 + 20t^7 = 1 - h_{01}(t)$$
$$h_{10}(t) = t - 20t^4 + 45t^5 - 36t^6 + 10t^7$$
$$h_{20}(t) = t^2 - 10t^4 + 20t^5 - 15t^6 + 4t^7$$
$$h_{30}(t) = t^3 - 4t^4 + 6t^5 - 4t^6 + t^7$$
$$h_{01}(t) = 35t^4 - 84t^5 + 70t^6 - 20t^7$$
$$h_{11}(t) = -15t^4 + 39t^5 - 34t^6 + 10t^7$$
$$h_{21}(t) = 5t^4 - 14t^5 + 13t^6 - 4t^7$$
$$h_{31}(t) = -t^4 + 3t^5 - 3t^6 + t^7$$

The derivative is

$$P'(x) = 140t^3(1 - 3t + 3t^2 - t^3)(y_{i+1} - y_i)/\Delta x + (1 - 80t^3 + 225t^4 - 216t^5 + 70t^6)v_i + t^3(-60 + 195t - 204t^2 + 70t^3)v_{i+1}$$
$$t\Delta x \left[(1 - 20t^2 + 50t^3 - 45t^4 + 14t^5)a_i + t^2(10 - 35t + 39t^2 - 14t^3)a_{i+1}\right]$$
$$\frac{t^2\Delta x^2}{6}\left[(3 - 16t + 30t^2 - 24t^3 + 7t^4)j_i + t(-4 + 15t - 18t^2 + 7t^3)j_{i+1}\right]$$

the second derivative is

$$P''(x) = 420t^2(1 - 4t + 5t^2 - 2t^3)(y_{i+1} - y_i)/\Delta x^2 + 60t^2(-4 + 15t - 18t^2 + 7t^3)v_i/\Delta x + 60t^2(-3 + 13t - 17t^2 + 7t^3)v_{i+1}/\Delta x$$
$$+ (1 - 60t^2 + 200t^3 - 225t^4 + 84t^5)a_i + t^2(30 - 140t + 195t^2 - 84t^3)a_{i+1}$$
$$+ t(1 - 8t + 20t^2 - 20t^3 + 7t^4)j_0\Delta x + t^2(-2 + 10t - 15t^2 + 7t^3)j_1\Delta x$$

Astute readers may have noticed that the Hermite basis requires more floating point operations (even Hornerized) than is strictly required for evaluation of a polynomial of the given degree. It is therefore tempting to convert the Hermite data to (say) Chebyshev coefficients, and use Reinch's modification to the Clenshaw recurrence[24] to evaluate the polynomial. This in fact can improve speed, but the conversion from the Hermite data to the Chebyshev coefficients is poorly conditioned. In our own case, this made Chebyshev coefficient determination in `long double` insufficient for evaluation in `double` precision, forcing use of quad precision. This made the constructor call unreasonably expensive, and the approach was abandoned. Note the incidental nature of this argument: there is still scope to use the Chebyshev basis if the user only requires computation of a single $p$ at a single precision, as precomputation and storage of the data in a `.cpp` is sensible in this case.